\documentclass[12pt]{article}

\usepackage{amsmath, epsfig, cite}
\usepackage{amssymb}
\usepackage{amsfonts}
\usepackage{latexsym}

\newtheorem{thm}{Theorem}[section]
\newtheorem{prop}[thm]{Proposition}

\newtheorem{conj}[thm]{Conjecture}
\newtheorem{lem}[thm]{Lemma}

\newcommand{\pf}{\noindent{\it Proof.} }

\numberwithin{equation}{section}

\newcommand{\qed}{{\hfill$\square$}\medskip}

\begin{document}

\nocite{*}
\begin{center}
{\Large\bf Proofs of two conjectures on Catalan triangle\\[5pt] numbers}
\end{center}

\vskip 2mm \centerline{Victor J. W. Guo and Xiuguo Lian}
\begin{center}
{School of Mathematical Sciences, Huaiyin Normal University, Huai'an 223300, Jiangsu,
 People's Republic of China\\
{\tt jwguo@hytc.edu.cn,   lianxiuguo@126.com} }

\end{center}


\vskip 0.7cm \noindent{\bf Abstract.} We prove two conjectures on sums of products of Catalan triangle numbers, which were originally conjectured by Miana, Ohtsuka, and Romero [Discrete Math. 340 (2017), 2388--2397].
The first one is proved by using Zeilberger's algorithm, and the second one is proved by establishing its $q$-analogue.


\vskip 2mm
\noindent{\it MR Subject Classifications}: 05A30; 05A10; 11B65.

\section{Introduction}
The well-known {\it Catalan numbers} $C_n:=\frac{1}{n+1}{2n\choose n}$ play an important role in combinatorics. For example, the Catalan number $C_n$
counts the number of sequences $a_1,a_2,\ldots, a_{2n}$ consisting of $n$ $1$'s and $n$ $(-1)$'s such that $a_1+\cdots+a_k\geqslant 0$ for $k=1,\ldots,2n$.
See \cite{Stanley1999,Stanley2015} for more combinatorial explanations of the Catalan numbers.

Shapiro \cite{Shapiro} introduced the {\it Catalan triangle} $(B_{n,k})_{n\geqslant k\geqslant 0}$ given by
$$
B_{n,k}:=\frac{k}{n}{2n\choose n-k}={2n-1\choose n-k}-{2n-1\choose n-k-1}.
$$
The Catalan triangle and related topics have been studied by many
authors over the past decade. See, for example, \cite{CC, Chu, GW,
GZ2010, GZ2011, GHJR, KP, MR2007, MR, MOR, SS, SM, ZP}. In
particular, Miana, Ohtsuka, and Romero \cite[Corollary 3.2(i),
Theorem 3.3]{MOR} proved that
\begin{align}
\sum_{k=0}^n B_{n,k}^3
&=\frac{1}{2}{2n\choose n}^3-\frac{3}{2}{2n\choose n}\sum_{k=n}^{2n-1}{k\choose n}{k\choose n-1}, \label{eq:recover}\\[5pt]
&=\frac{1}{2n}{2n\choose n}\sum_{k=1}^{n}k{2n-k-1\choose n-1}^2. \notag
\end{align}
They \cite[Conjecture 4.2]{MOR} made the following conjecture:
\begin{conj}\label{conj:one}
For all positive integers $m$ and $n$, the identity
\begin{align}
\sum_{k=0}^{r}B_{n,k}^2 B_{m,k}=\frac{1}{2}{2n\choose n}^2 {2m\choose m}
\left[ 1-\frac{n+2m}{r}{n+m\choose n}^{-1}{n+r\choose n}^{-1}\sum_{k=0}^{r-1}{s+k\choose s}{n+k\choose n-1}\right]  \label{eq:one}
\end{align}
holds for $r=\min\{n,m\}$ and $s=\max\{n,m\}$. In the particular case $m=n$, we recover \eqref{eq:recover}.
\end{conj}

In this paper, we shall prove the following stronger version of Conjecture \ref{conj:one}.
\begin{thm}\label{thm:first}
The identity \eqref{eq:one} holds for $r=m$ and $s=n$. It also holds for $r=n$ and $s=m$.
\end{thm}

We shall also give the following new identity.
\begin{thm}\label{thm:second}
Let $m$ and $n$ be non-negative integers. Then
\begin{align}
\sum_{k=1}^{m}{n+k-2\choose k-1}{n-1\choose k-1}{m+n\choose m-k}
=\sum_{k=0}^m k{m+n-k-1\choose n-1}^2. \label{eq:new-id}
\end{align}
\end{thm}

Miana, Ohtsuka, and Romero \cite{MOR} considered the numbers
$$
C_{n,k}:=\frac{n-2k}{n}{n\choose k}={n-1\choose k}-{n-1\choose k-1},
$$
and called them {\it Catalan triangle numbers}. It is easy to see that $B_{n,k}=C_{2n,n-k}$. Miana, Ohtsuka, and Romero \cite[ Theorem 1.2(i) and Theorem 3.1(i)]{MOR} proved that
\begin{align*}
\sum_{k=0}^{a}C_{n,k}&={n-1\choose a},\\[5pt]
\sum_{k=0}^{a}C_{n,k}^3&=4{n-1\choose a}^3-3{n-1\choose a}\sum_{j=0}^{n-1}{j\choose a}{j\choose n-a-1},
\end{align*}
and proposed the following conjecture \cite[Conjecture 4.1]{MOR}:
\begin{conj}\label{conj:2}
Let $n$ and $a$ be positive integers with $n>a$, and let $r$ be a non-negative integer. Then
\begin{align}
\sum_{k=0}^{a}C_{n,k}^{2r+1}\equiv 0\pmod{{n-1\choose a}}.  \label{eq:cnk}
\end{align}
\end{conj}

Note that Conjecture \ref{conj:2} is true for some special cases (see Guo and Zeng \cite[Theorem 1.4]{GZ2010}, Guo and Zeng \cite[Theorem 1.4]{GZ2011}, and Guo and Wang \cite[Theorem 1.3]{GW}).

It is easy to see that $C_{n,k}=-C_{n,n-k}$. Hence, to prove \eqref{eq:cnk}, it suffices to prove \eqref{eq:cnk} for the case where $n\geqslant 2a$.
It is easy to see that the congruence \eqref{eq:cnk} can be written as the following two congruences according to the parity of $n$ (also reversing the summation order):
\begin{align}
\sum_{k=a}^{n}B_{n,k}^{2r+1}\equiv 0\pmod{{2n-1\choose n-a}},\quad 0\leqslant a\leqslant n,  \label{eq:bnk}\\[5pt]
\sum_{k=a}^{n}A_{n,k}^{2r+1}\equiv 0\pmod{{2n\choose n-a}},\quad 0\leqslant a\leqslant n,  \label{eq:ank}
\end{align}
where $$A_{n,k}:=\frac{2k+1}{2n+1}{2n+1\choose n-k},\quad 0\leqslant k\leqslant n.$$

In this paper, we shall confirm Conjecture \ref{conj:2} by establishing \eqref{eq:bnk} and \eqref{eq:ank}, respectively.
More precisely, we shall give $q$-analogues of \eqref{eq:bnk} and \eqref{eq:ank}.

The {\it $q$-integers} are defined by $[n]=\frac{1-q^n}{1-q}$ and the {\it $q$-shifted factorials} (see \cite{GR}) are defined by $(a;q)_0=1$
and $(a;q)_n=(1-a)(1-aq)\cdots (1-aq^{n-1})$ for $n=1,2,\ldots.$ The {\it $q$-binomial coefficients} are defined as
$$
{n\brack k}={n\brack k}_q=
\begin{cases}\displaystyle\frac{(q;q)_n}{(q;q)_k (q;q)_{n-k}}, &\text{if $0\leqslant k\leqslant n$,} \\[10pt]
0, &\text{otherwise.}
\end{cases}
$$
It is well-known that the $q$-binomial coefficients are polynomials in $q$ with integer coefficients. The polynomials $B_{n,k}(q)$ and $A_{n,k}(q)$ are given by (see \cite{GZ2010,GW2})
\begin{align}
B_{n,k}(q)&:=\frac{[k]}{[n]}{2n\brack n-k},\quad 1\leqslant k\leqslant n, \notag\\[5pt]
A_{n,k}(q)&:=q^{n-k}\frac{[2k+1]}{[2n+1]}{2n+1\brack n-k}={2n\brack n-k}-{2n\brack n-k-1},\quad 0\leqslant k\leqslant n.  \label{eeq:ank-integer}
\end{align}

Let $P(q)$ be a Laurent polynomial in $q$ and let $D(q)$ be a polynomial in $q$. For convenience, we write
$$
P(q)\equiv 0\pmod{D(q)}
$$
if $P(q)/D(q)$ is still a Laurent polynomial in $q$.

We have the following $q$-versions of \eqref{eq:bnk} and \eqref{eq:ank}.
\begin{thm}\label{thm:bnk-power}
Let $n$ be a positive integer,  and let $a$ and $r$ be non-negative
integers with $a\leqslant n$. Then, for $0\leqslant j\leqslant
2r+1$, there holds
\begin{align*}
\sum_{k=a}^{n}(1+q^k)B_{n,k}(q)^{2r+1} q^{jk^2-(r+1)k} \equiv 0\pmod{(1+q^n){2n-1\brack n-a}}.
\end{align*}
\end{thm}

\begin{thm}\label{thm:ank-power}
Let $n$ be a positive integer,  and let $a$ and $r$ be non-negative integers with $a\leqslant n$. Then, for $0\leqslant j\leqslant 2r+1$, there holds
\begin{align*}
\sum_{k=a}^{n} A_{n,k}(q)^{2r+1} q^{j(k^2+k)}  \equiv 0\pmod{{2n\brack n-a}}.
\end{align*}
\end{thm}

It is clear that \eqref{eq:bnk} and \eqref{eq:ank} follow from Theorems \ref{thm:bnk-power} and \ref{thm:ank-power} by letting $q=1$. It should also
be mentioned that the $a=0$ case of  Theorem \ref{thm:bnk-power} follows from \cite[Theorem 1.3]{GW}, while the $a=0$ case of  Theorem \ref{thm:ank-power} can be deduced from \cite[Theorem 1.3]{GW2}.

\section{Proof of Theorem \ref{thm:first}}
We first prove \eqref{eq:one} holds for $r=m$ and $s=n$. Namely,
\begin{align}
\sum_{k=0}^{m}B_{n,k}^2 B_{m,k}=\frac{1}{2}{2n\choose n}^2 {2m\choose m}
\left[ 1-\frac{n+2m}{m}{n+m\choose n}^{-2}\sum_{k=0}^{m-1}{n+k\choose n}{n+k\choose n-1}\right].  \label{eq:one-02}
\end{align}
It can be easily proved by induction that (see \cite[Theorem 2.3(i)]{MOR})
$$
{n+m\choose n}^2=\sum_{k=0}^{m}\frac{2m+n-2k}{n}{m+n-k-1\choose n-1}^2.
$$
Moreover, by reversing the summation order, we have
$$
\sum_{k=0}^{m-1}{n+k\choose n}{n+k\choose n-1}=\sum_{k=0}^{m}\frac{m-k}{n}{m+n-k-1\choose n-1}^2.
$$
Hence, the identity  \eqref{eq:one-02} is equivalent to
\begin{align}
2m{2n\choose n}^{-2} {2m\choose m}^{-1}{n+m\choose n}^{2}\sum_{k=0}^{m}B_{n,k}^2 B_{m,k}
=\sum_{k=0}^{m-1}k{m+n-k-1\choose n-1}^2. \label{eq:one-2}
\end{align}
Denote the left-hand side and right-hand side of \eqref{eq:one-2} by $S_n(m)$ and $T_n(m)$, respectively. Then Zeilberger's algorithm \cite{Koepf,PWZ} gives
\begin{align*}
(2m+n)S_n(m+1)&=(2m+n+2)S_n(m)+n{n+m\choose n}^2,\\[5pt]
(2m+n)T_n(m+1)&=(2m+n+2)T_n(m)+n{n+m\choose n}^2.
\end{align*}
Since $S_n(0)=T_n(0)=0$, we conclude that $S_n(m)=T_n(m)$ for all non-negative integers $m$.

To prove that  \eqref{eq:one} also holds for $r=n$ and $s=m$, noting that
$$
\sum_{k=0}^{m}B_{n,k}^2 B_{m,k}=\sum_{k=0}^{n}B_{n,k}^2 B_{m,k},
$$
it suffices to prove
\begin{align}
\frac{1}{m}{n+m\choose n}^{-1}\sum_{k=0}^{m-1}{n+k\choose n}{n+k\choose n-1}
=\frac{1}{n}{2n\choose n}^{-1}\sum_{k=0}^{n-1}{m+k\choose m}{n+k\choose n-1}, \label{eq:one-suff}
\end{align}
which is equivalent to
\begin{align*}
&\frac{m+1}{n}{n+m+1\choose n}{2n\choose n}^{-1}\sum_{k=0}^{n-1}{m+k+1\choose m+1}{n+k\choose n-1} \\[5pt]
&{}-\frac{m}{n}{n+m\choose n}{2n\choose n}^{-1}\sum_{k=0}^{n-1}{m+k\choose m}{n+k\choose n-1} ={n+m\choose n}{n+m\choose n-1},
\end{align*}
or
\begin{align}
\sum_{k=0}^{n-1}\frac{(n+k+1)m+(n+1)(k+1)}{m+1}{m+k\choose m}{n+k\choose n-1}=n{n+m\choose n-1}{2n\choose n}. \label{eq:one-3}
\end{align}
But \eqref{eq:one-3} can be easily proved by Gosper's algorithm. This completes the proof of \eqref{eq:one-suff}.

\section{Proof of Theorem \ref{thm:second}}
The $m=3$ case of \cite[Corollary 4.1]{GZ2010} gives
\begin{align*}
&\sum_{k=1}^{n_1}k^3 {n_1+n_2\choose n_1+k}{n_2+n_3\choose n_2+k}{n_3+n_1\choose n_3+k} \\[5pt]
&\quad =\frac{n_1 n_3}{2}{n_1+n_3\choose n_1}\sum_{k=1}^{n_1}{n_3+k-2\choose k-1}{n_1-1\choose k-1}{n_2+n_3\choose n_2-k}.
\end{align*}
Namely,
\begin{align}
&\sum_{k=1}^{n_1}\frac{k^3}{n_1n_2n_3} {2n_1\choose n_1+k}{2n_2\choose n_2+k}{2n_3\choose n_3+k} \notag \\[5pt]
&\quad =\frac{1}{2n_2}{n_1+n_3\choose n_1}\frac{(n_1+n_2)!(n_2+n_3)!(n_3+n_1)!}{(2n_1)!(2n_2)!(2n_3)!}\sum_{k=1}^{n_1}{n_3+k-2\choose k-1}{n_1-1\choose k-1}{n_2+n_3\choose n_2-k}. \label{eq:n123}
\end{align}
Letting $n_1=n_3=n$ and $n_2=m$ in \eqref{eq:n123}, we obtain
\begin{align}
2m{2n\choose n}^{-2} {2m\choose m}^{-1}{n+m\choose n}^{2}\sum_{k=0}^{m}B_{n,k}^2 B_{m,k}
=\sum_{k=1}^{n}{n+k-2\choose k-1}{n-1\choose k-1}{m+n\choose m-k}. \label{eq:bbb-spec}
\end{align}
Combining \eqref{eq:one-2} and \eqref{eq:bbb-spec}, we get \eqref{eq:new-id}.

\medskip
\noindent{\it Remark.} The identity \eqref{eq:new-id} can also be proved by Zeilberger's algorithm.

\section{Proof of Theorem \ref{thm:bnk-power}}

We first need to establish the following result, which is similar to \cite[Theorem 1.1]{GW}.
\begin{thm}\label{thm:half}
Let $n$ be a positive integer and let $a,r$ be non-negative integers with $a\leqslant n$. Then, for $j=0,1$, there holds
\begin{align*}
\sum_{k=a}^{n}[2k][k]^{2r}q^{(r+1)(n-k)+jk^2}{2n\brack n-k} \equiv 0\pmod{[n+a]{2n\brack n-a}}.
\end{align*}
\end{thm}
\pf
For $j=0$, let
\begin{align*}
S_r(a,n;q)=\sum_{k=a}^{n}[2k][k]^{2r}q^{(r+1)(n-k)}{2n\brack n-k}.
\end{align*}
It is easily seen that (via telescoping)
\begin{align*}
S_0(a,n;q)=\sum_{k=a}^n[2k]q^{n-k}{2n\brack n-k}
&=\sum_{k=a}^n [n-k+1]{2n\brack n-k+1}
-\sum_{k=a}^n [n-k]{2n\brack n-k}\\
&=[n-a+1]{2n\brack n-a+1}\\
&=[n+a]{2n\brack n-a}.
\end{align*}

For $r\geqslant 1$, since
\begin{align*}
[k]^2{2n\brack n-k}q^{n-k}=[n]^2{2n\brack n-k}-[2n][2n-1]{2n-2\brack n-k-1},
\end{align*}
we have
\begin{align}
S_r(a,n;q)=[n]^2S_{r-1}(a,n;q)-q^r[2n][2n-1] S_{r-1}(a,n-1;q),\ n=1,2,\ldots. \label{eq:sumnkl}
\end{align}
Applying the recurrence relation \eqref{eq:sumnkl} and by induction on $r$,
we can easily prove that, for all positive integers $r$, there holds
\begin{align}
S_r(a,n;q)\equiv 0\pmod{[n+a]{2n\brack n-a}}. \label{eq:srnq}
\end{align}

For $j=1$, let
\begin{align*}
T_r(a,n;q)=\sum_{k=a}^{n}[2k][k]^{2r}q^{(r+1)(n-k)+k^2}{2n\brack n-k}.
\end{align*}
Using the relation ${n\brack k}_{q^{-1}}={n\brack k}q^{k(k-n)}$, we have
$ T_r(a,n;q)=q^{n^2+2rn+2n-2r-1} S_r(a,n;q^{-1})$. Therefore, from \eqref{eq:srnq} we deduce that
\begin{align*}
T_r(a,n;q)\equiv 0\pmod{[n+a]{2n\brack n-a}}.
\end{align*}
Since both $S_r(a,n;q)$ and $T_r(a,n;q)$ are polynomials in $q$, we obtain the desired result.
\qed

We also need the following generalization of Theorem \ref{thm:half}.
\begin{thm}\label{thm:nnnhalf}
Let $n_1,\ldots,n_{m},n_{m+1}=n_1$ be positive integers. Then for any non-negative integers $a$, $j$ and $r$ with $a\leqslant n_1$ and $j\leqslant m$, the expression
\begin{align}
\frac{1}{[n_m+a]}{n_1+n_{m}\brack n_1-a}^{-1}
\sum_{k=a}^{n_1}[2k][k]^{2r}q^{jk^2-(r+1)k}\prod_{i=1}^{m} {n_i+n_{i+1}\brack n_i+k}   \label{eq:main-two}
\end{align}
is a Laurent polynomial in $q$ with integer coefficients.
\end{thm}
\pf
Denote \eqref{eq:main-two} by $S_r(a;n_1,\ldots,n_{m};j,q)$. Let
$$
C(a_1,\ldots,a_l;k)=\prod_{i=1}^l {a_i+a_{i+1}\brack a_i+k}\quad (a_{l+1}=a_1).
$$
Then we can write $S_a(n_1,\ldots,n_{m};r,j,q)$ as
\begin{equation}\label{eq:rewriting}
S_r(a;n_1,\ldots,n_{m};j,q)
=\frac{(q;q)_{n_1-a} (q;q)_{n_{m}+a-1}}{(q;q)_{n_1+n_{m}}}
\sum_{k=a}^{n_1} [2k][k]^{2r}q^{jk^2-(r+1)k} C(n_1,\ldots,n_{m};k).
\end{equation}

For $m\geqslant 3$, there holds
\begin{align}
C(n_1,\ldots,n_m;k)=\frac{(q;q)_{n_2+n_3}(q;q)_{n_m+n_1}}{(q;q)_{n_1+n_2}(q;q)_{n_m+n_3}}
{n_1+n_2\brack n_1+k}{n_1+n_2\brack n_2+k}C(n_3,\ldots,n_m;k).  \label{eq:cn1nk}
\end{align}
By the $q$-Chu-Vandermonde summation formula (see, for example, \cite[p.~37, (3.3.10)]{Andrews}), we have
$$
{n_1+n_2\brack n_1+k}=\sum_{s=0}^{n_1-k}{n_1-k\brack s}{n_2+k\brack s+2k} q^{s(s+2k)},
$$
which may be rewritten as
\begin{align}
{n_1+n_2\brack n_1+k}{n_1+n_2\brack n_2+k}
=\sum_{s=0}^{n_1-k}\frac{q^{s^2+2ks}(q;q)_{n_1+n_2}}{(q;q)_{s}(q;q)_{s+2k}(q;q)_{n_1-k-s}(q;q)_{n_2-k-s}},\label{eq:qpfaffn1n2}
\end{align}
where we assume that $\frac{1}{(q;q)_n}=0$ for any negative integer $n$.
Substituting \eqref{eq:cn1nk} and \eqref{eq:qpfaffn1n2} into \eqref{eq:rewriting}, we obtain
\begin{align*}
&\hskip -2mm S_r(a;n_1,\ldots,n_{m};j,q) \\
&=\frac{(q;q)_{n_2+n_3}(q;q)_{n_1-a} (q;q)_{n_m+a-1}}{(q;q)_{n_m+n_3}}\sum_{k=a}^{n_1}\sum_{s=0}^{n_1-k}
\frac{q^{s^2+2ks+jk^2-(r+1)k}[2k][k]^{2r} C(n_3,\ldots,n_{m};k)}{(q;q)_{s}(q;q)_{s+2k}(q;q)_{n_1-k-s}(q;q)_{n_2-k-s}}  \\
&=\frac{(q;q)_{n_2+n_3}(q;q)_{n_1-a} (q;q)_{n_m+a-1}}{(q;q)_{n_m+n_3}}\sum_{l=a}^{n_1} q^{l^2}\sum_{k=a}^{l}
\frac{q^{(j-1)k^2-(r+1)k}[2k][k]^{2r} C(n_3,\ldots,n_{m};k)}{(q;q)_{l-k}(q;q)_{l+k}(q;q)_{n_1-l}(q;q)_{n_2-l}},
\end{align*}
where $l=s+k$. Noticing that, for $m\geqslant 3$,
$$
\frac{C(n_3,\ldots, n_{m};k)}{(q;q)_{l-k}(q;q)_{l+k}}
=\frac{(q;q)_{n_{m}+n_3}}{(q;q)_{n_3+l}(q;q)_{n_{m}+l}}C(l,n_3,\ldots, n_{m};k),
$$
we are led to the following recurrence relation
\begin{align}
S_r(a;n_1,\ldots,n_{m};j,q)
=\sum_{l=a}^{n_1} q^{l^2}{n_1-a\brack l-a}{n_2+n_3\brack n_2-l} S_r(a;l,n_3,\ldots,n_{m};j-1,q). \label{eq:recsr}
\end{align}

Similarly, for $m=2$, we have
\begin{align}
S_r(a;n_1,n_{2};j,q)=\sum_{l=a}^{n_1} q^{l^2}{n_1-a\brack l-a}{n_2+a-1\brack l+a-1}S_r(a;l;j-1,q).  \label{eq:sn1n2}
\end{align}

Now we can give an inductive proof of the theorem.  For $m=1$, the conclusion is true by Theorem \ref{thm:half}.
Suppose that $S_r(a;n_1,\ldots,n_{m-1};j,q)$ is a Laurent polynomial in $q$ with integer coefficients for some $m\geqslant 2$ and all $j$ with $0\leqslant j\leqslant m-1$.
Then by \eqref{eq:recsr} (if $m\geqslant 3$) or \eqref{eq:sn1n2} (if $m=2$), so is $S_r(a;n_1,\ldots,n_{m};j,q)$ for  $1\leqslant j\leqslant m$. Furthermore, since
\begin{align*}
S_r(a;n_1,\ldots,n_{m};0,q)=S_r(a;n_1,\ldots,n_{m};m,q^{-1}) q^{n_1n_2+n_2n_3+\cdots+n_{m-1}n_m-n_m+a(a+n_m-n_1-1)-2r}
\end{align*}
for $m\geqslant 2$, one sees  that $S_r(a;n_1,\ldots,n_{m};0,q)$ is also a Laurent polynomial in $q$ with integer coefficients. This completes the proof.
\qed

Note that, by \eqref{eq:sn1n2} and the $q$-Chu-Vandermonde summation formula, there holds
\begin{align*}
S_0(a;n_1,n_{2};1,q)
={n_1+n_2-1\brack n_1+a-1}q^{a^2-a}.
\end{align*}

Let $\Phi_n(x)$ be the $n$-th {\it cyclotomic polynomial}. The following result is very useful in dealing with $q$-binomial coefficients (see \cite[(10)]{KW} or \cite{CH,GZ06}).
\begin{prop}\label{prop:factor}
The $q$-binomial coefficient ${m\brack k}$ can be factorized into
$$
{m\brack k}=\prod_{d}\Phi_d(q),
$$
where the product is over all positive integers $d\leqslant m$ such that
$\lfloor k/d\rfloor+\lfloor (m-k)/d\rfloor<\lfloor m/d\rfloor$.
\end{prop}

\medskip
\noindent{\it Proof of Theorem {\rm\ref{thm:bnk-power}}.}
Letting $m=2r+1$ and $n_1=\cdots=n_{2r+1}=n$ in Theorem~\ref{thm:nnnhalf}, we see that
\begin{align*}
&\frac{1}{[n+a]}{2n\brack n-a}^{-1}
\sum_{k=a}^{n}[2k][k]^{2r}q^{jk^2-(r+1)k} {2n\brack n-k}^{2r+1} \\[5pt]
&\quad=\frac{1}{[2n]}{2n-1\brack n-a}^{-1}
\sum_{k=a}^{n}(1+q^k)q^{jk^2-(r+1)k} \left([n]B_{n,k}(q)\right)^{2r+1}
\end{align*}
is a Laurent polynomial in $q$ with integer coefficients.  Note that $B_{n,k}(q)$ is a polynomial in $q$ with integer coefficients (see \cite{GZ2010,GK}). Therefore,
\begin{align*}
\frac{\gcd\left([2n]{2n-1\brack n-a},[n]^{2r+1}\right)}{[2n]{2n-1\brack n-a}}\sum_{k=a}^{n}(1+q^k)q^{jk^2-(r+1)k} B_{n,k}(q)^{2r+1}
\end{align*}
is a Laurent polynomial in $q$ with integer coefficients.

It is well-known that
\begin{align}
[n]=\prod_{\substack{d|n\\ d>1}}\Phi_d(q),  \label{eq:qint-factor}
\end{align}
Therefore, from Proposition \ref{prop:factor} we immediately deduce that
\begin{align*}
\gcd\left({2n-1\brack n-a},[n]\right)=1.  
\end{align*}
Moreover, we have $[2n]=(1+q^n)[n]$ and $\gcd(1+q^n,[n])=1$. It follows that
$$
\frac{\gcd\left([2n]{2n-1\brack n-a},[n]^{2r+1}\right)}{[2n]{2n-1\brack n-a}}=\frac{1}{(1+q^n){2n-1\brack n-a}}.
$$
This completes the proof.
\qed

\section{Proof of Theorem \ref{thm:ank-power}}
We first need to establish the following result.
\begin{lem}Let $n$ be a positive integer and let $a,s$ be non-negative integers with $a\leqslant n$. Then
\begin{align}
\sum_{k=a}^{n}q^{n-k}[2k+1] {2n+1\brack n-k}(q^{-k};q)_s (q^{k+1};q)_s \equiv 0\pmod{[2n+1]{2n\brack n-a}}.\label{eq:r1m1-1}
\end{align}
\end{lem}
\pf We proceed by induction on $s$. For $s=0$, we have
\begin{align*}
\sum_{k=a}^{n}q^{n-k}[2k+1] {2n+1\brack n-k}
&=[2n+1]\sum_{k=a}^{n}\left({2n\brack n-k}-{2n\brack n-k-1}\right) \\[5pt]
&=[2n+1]{2n\brack n-a}.
\end{align*}
Suppose \eqref{eq:r1m1-1} is true for $s$. Noticing the relations
\begin{align*}
&\hskip -2mm {2n+1\brack n-k}(q^{-k};q)_{s+1} (q^{k+1};q)_{s+1}  \\
&=(1-q^{s-n})(1-q^{s+n+1}){2n+1\brack n-k}(q^{-k};q)_{s} (q^{k+1};q)_{s}  \\
&\quad{}+q^{s-n}(1-q^{2n})(1-q^{2n+1}){2n-1\brack n-k-1}(q^{-k};q)_{s} (q^{k+1};q)_{s}
\end{align*}
and
\begin{align}
[2n][2n+1][2n-1]{2n-2\brack n-a-1}=[2n+1]{2n\brack n-a}[n-a][n+a],  \label{eq:relation-2}
\end{align}
we can easily prove that \eqref{eq:r1m1-1} holds for $s+1$.
\qed

We have the following generalization of \eqref{eq:r1m1-1}.
\begin{lem}
Let $n$ be a positive integer and let $a,r,s$ be non-negative integers with $a\leqslant n$. Then
\begin{align}
&\sum_{k=a}^{n}q^{-(2r+1)k}[2k+1]^{2r+1} {2n+1\brack n-k}(q^{-k};q)_s (q^{k+1};q)_s \equiv 0\pmod{[2n+1]{2n\brack n-a}}. \label{eq:mod-r1m1-1}
\end{align}
\end{lem}
\pf We proceed by induction on $r$.
Denote the left-hand side of \eqref{eq:mod-r1m1-1} by $X_{r}(a,n,s;q)$.
By \eqref{eq:r1m1-1}, one sees that \eqref{eq:mod-r1m1-1} is true for $r=0$.
For $r\geqslant 1$, suppose that
$$
X_{r-1}(a,n,s;q)\equiv 0\mod [2n+1]{2n\brack n-a}
$$
holds for any positive integer $n$ and non-negative integers $a$, $s$ with $a\leqslant n$.
It is easy to check that
\begin{align*}
{2n+1\brack n-k}[2k+1]^2
&=q^{2k-2n}{2n+1\brack n-k}[2n+1]^2 \\[5pt]
&\quad{}-q^{2k-2n}{2n-1\brack n-k-1}[2n][2n+1](1+q^{n-s})(1+q^{n+s+1}) \\[5pt]
&\quad{}+q^{2k-n-s}{2n-1\brack n-k-1}[2n][2n+1](1-q^{s-k})(1-q^{s+k+1}),
\end{align*}
and therefore,
\begin{align}
X_{r}(a,n,s;q)
&=q^{-2n}[2n+1]^2X_{r-1}(a,n,s;q)\notag\\[5pt]
&\quad{}-q^{-2n}[2n][2n+1](1+q^{n-s})(1+q^{n+s+1})X_{r-1}(a,n-1,s;q) \notag\\[5pt]
&\quad{}+q^{-n-s}[2n][2n+1]X_{r-1}(a,n-1,s+1;q).  \label{eq:rec-3-term}
\end{align}
By the induction hypothesis and applying \eqref{eq:relation-2},
we immediately deduce from the recurrence \eqref{eq:rec-3-term} that
\eqref{eq:mod-r1m1-1} holds for $r$.   \qed

\begin{thm}\label{thm:fatorodd-second}
Let $n_1,\ldots,n_{m},n_{m+1}=n_1$ be positive integers. Then for any non-negative integers $a$, $j$ and $r$ with $a\leqslant n_1$ and $j\leqslant m$, the expression
\begin{align}
\frac{1}{[n_1+n_m+1]}{n_1+n_{m}\brack n_1-a}^{-1}\sum_{k=a}^{n_1} q^{j(k^2+k)-(2r+1)k}[2k+1]^{2r+1}\prod_{i=1}^m {n_i+n_{i+1}+1\brack n_i-k} \label{eq:ank-q}
\end{align}
is a Laurent polynomial in $q$ with integer coefficients.
\end{thm}
\pf Denote \eqref{eq:ank-q} by $\overline{S}_r(a;n_1,\ldots,n_{m};j,q)$.  Let
$$
\overline{C}(a_1,\ldots,a_l;k)=\prod_{i=1}^l {a_i+a_{i+1}+1\brack a_i-k},
$$
where $a_{l+1}=a_1$. Then
\begin{align}
&\hskip -2mm \overline{S}_r(a;n_1,\ldots,n_{m};j,q)  \notag \\[5pt]
&=\frac{(q;q)_{n_1-a}(q;q)_{n_m+a}}{(q;q)_{n_1+n_m+1}}
\sum_{k=a}^{n_1} q^{j(k^2+k)-(2r+1)k}[2k+1]^{2r+1} \overline{C}(n_1,\ldots,n_{m};k),  \label{eq:sr-n1nm}
\end{align}

For $m\geqslant 3$, we have
\begin{align}
\overline{C}(n_1,\ldots,n_m;k)=\frac{(q;q)_{n_2+n_3+1}(q;q)_{n_m+n_1+1}}{(q;q)_{n_1+k+1}(q;q)_{n_2-k}(q;q)_{n_m+n_3+1}}
{n_1+n_2+1\brack n_1-k}\overline{C}(n_3,\ldots,n_m;k).  \label{eq:C}
\end{align}
Applying \eqref{eq:C} and the $q$-Chu-Vandermonde summation formula (see \cite[p.~37, (3.3.10)]{Andrews})
\begin{align}
{n_1+n_2+1\brack n_1-k}
=\sum_{s=0}^{n_1-k}\frac{q^{s(s+2k+1)}(q;q)_{n_1+k+1}(q;q)_{n_2-k}}{(q;q)_s (q;q)_{s+2k+1}(q;q)_{n_1-k-s}(q;q)_{n_2-k-s}},  \label{eq:vandermonde}
\end{align}
we may write \eqref{eq:sr-n1nm} as
\begin{align*}
&\overline{S}_r(a;n_1,\ldots,n_{m};j,q) \\
&=\frac{(q;q)_{n_2+n_3+1}(q;q)_{n_1-a} (q;q)_{n_m+a}}{(q;q)_{n_m+n_3+1}}\\[5pt]
&\quad\times\sum_{k=a}^{n_1}\sum_{s=0}^{n_1-k}
\frac{q^{j(k^2+k)-(2r+1)k+s(s+2k+1)}[2k+1]^{2r+1} \overline{C}(n_3,\ldots,n_{m};k)}{(q;q)_{s}(q;q)_{s+2k+1}(q;q)_{n_1-k-s}(q;q)_{n_2-k-s}}  \\
&=\frac{(q;q)_{n_2+n_3+1}(q;q)_{n_1-a} (q;q)_{n_m+a}}{(q;q)_{n_m+n_3+1}}\\[5pt]
&\quad\times\sum_{l=a}^{n_1} q^{l^2+l}\sum_{k=a}^{l}
\frac{q^{(j-1)(k^2+k)-(2r+1)k}[2k+1]^{2r+1} \overline{C}(n_3,\ldots,n_{m};k)}{(q;q)_{l-k}(q;q)_{l+k}(q;q)_{n_1-l}(q;q)_{n_2-l}},
\end{align*}
where $l=s+k$. Noticing that
$$
\frac{\overline{C}(n_3,\ldots, n_{m};k)}{(q;q)_{l-k}(q;q)_{l+k+1}}
=\frac{(q;q)_{n_{m}+n_3+1}}{(q;q)_{n_3+l+1}(q;q)_{n_{m}+l+1}}\overline{C}(l,n_3,\ldots, n_{m};k),
$$
we obtain
\begin{align}
\overline{S}_r(a;n_1,\ldots,n_{m};j,q)
=\sum_{l=a}^{n_1} q^{l^2+l}{n_1-a\brack l-a}{n_2+n_3+1\brack n_2-l} \overline{S}_r(a;l,n_3,\ldots,n_{m};j-1,q) \label{eq:S-recsr}
\end{align}
for $m\geqslant 3$. Similarly, for $m=2$, applying \eqref{eq:vandermonde} we get
\begin{align}
\overline{S}_r(a;n_1,n_{2};j,q)=\sum_{l=a}^{n_1} q^{l^2+l}{n_1-a\brack l-a}{n_2+a\brack l+a}\overline{S}_r(a;l;j-1,q). \label{eq:Sn1n2}
\end{align}

Similarly to the inductive proof of Theorem \ref{thm:nnnhalf}, using \eqref{eq:mod-r1m1-1} (with $s=0$), \eqref{eq:S-recsr},  \eqref{eq:Sn1n2} and the following relation
\begin{align*}
\overline{S}_r(a;n_1,\ldots,n_{m};0,q)=\overline{S}_r(a;n_1,\ldots,n_{m};m,q^{-1}) q^{n_1n_2+n_2n_3+\cdots+n_{m-1}n_m-n_1-n_m+a(a+n_m-n_1)}
\end{align*}
for $m\geqslant 2$, we can prove the theorem.  \qed

Note that, by \eqref{eq:Sn1n2} and the $q$-Chu-Vandermonde summation formula, there holds
\begin{align*}
\overline{S}_0(a;n_1,n_{2};1,q)
={n_1+n_2\brack n_1+a}q^{a^2}.
\end{align*}

\medskip
\noindent{\it Proof of Theorem {\rm\ref{thm:ank-power}}.}
Letting $m=2r+1$ and $n_1=\cdots=n_{2r+1}=n$ in Theorem \ref{thm:fatorodd-second}, we see that
\begin{align*}
&\hskip -2mm\frac{1}{[2n+1]}{2n\brack n-a}^{-1}\sum_{k=a}^{n} q^{j(k^2+k)-(2r+1)k}[2k+1]^{2r+1}{2n+1\brack n-k}^{2r+1} \\[5pt]
&=\frac{1}{[2n+1]}{2n\brack n-a}^{-1}\sum_{k=a}^{n} q^{j(k^2+k)-(2r+1)n}[2n+1]^{2r+1}A_{n,k}(q)^{2r+1}
\end{align*}
is a Laurent polynomial in $q$ with integer coefficients. By \eqref{eeq:ank-integer}, we know that $A_{n,k}(q)$ is a polynomial in $q$ with integer coefficients.
It follows that
\begin{align*}
\frac{\gcd\left([2n+1]{2n\brack n-a},[2n+1]^{2r+1}\right)}{[2n+1]{2n\brack n-a}}\sum_{k=a}^{n} A_{n,k}(q)^{2r+1} q^{j(k^2+k)}
\end{align*}
is a Laurent polynomial in $q$ with integer coefficients. By \eqref{eq:qint-factor} and Proposition \ref{prop:factor}, it is easy to see that
\begin{align*}
\gcd\left({2n\brack n-a},[2n+1]\right)=1, 
\end{align*}
and so
$$
\frac{\gcd\left([2n+1]{2n\brack n-a},[2n+1]^{2r+1}\right)}{[2n+1]{2n\brack n-a}}=\frac{1}{{2n\brack n-a}}.
$$
This completes the proof. \qed

\section{Two open problems}
In this section, we give two related conjectures for further study. The first one is a stronger version of Theorem \ref{thm:nnnhalf}, and is also a generalization of \cite[Conjecture 6.3]{GW}.
\begin{conj}
Let $n_1,\ldots,n_{m},n_{m+1}=n_1$ be positive integers. Then, for any integer $j$ and  non-negative integers $a,r$ with $a\leqslant n_1$, the expression
\begin{align*}
\frac{1}{[n_m+a]}{n_1+n_{m}\brack n_1-a}^{-1}
\sum_{k=a}^{n_1}[2k][k]^{2r}q^{jk^2-(r+1)k}\prod_{i=1}^{m} {n_i+n_{i+1}\brack n_i+k}
\end{align*}
is a Laurent polynomial in $q$, and it has non-negative integer coefficients if $0\leqslant j\leqslant m$.
\end{conj}

\begin{conj}
Theorem {\rm\ref{thm:fatorodd-second}} is still true for any integer $j$.
\end{conj}

\noindent{\bf Acknowledgments.} The authors thank the anonymous referee for a careful reading of this paper.
The first author was partially supported by the National Natural Science Foundation of China (grant 11771175),
the Natural Science Foundation of Jiangsu Province (grant BK20161304),
and the Qing Lan Project of Education Committee of Jiangsu Province.

\end{document}